\documentclass[12pt,a4paper]{amsart}
\usepackage[cp1250]{inputenc}
\usepackage{amssymb}

\theoremstyle{plain}
\newtheorem*{tw}{Theorem}
\newtheorem{wn}{Corollary}
\newtheorem{df}{Definition}
\newtheorem{uw}{Remark}
\newtheorem{fa}{Fact}

\newcommand{\CC}{\mathbb{C}}

\newcommand{\ie}{i.\,e. }
\newcommand{\floor}[1]{\lfloor{#1}\rfloor}
\DeclareMathOperator{\NWW}{lcm}
\DeclareMathOperator{\NWD}{gcd}
\DeclareMathOperator{\mdeg}{mdeg}

\title{On multidegrees of polynomial automorphisms of $\mathbb{C}^3$}
\author[J. Zygad{\l}o]{Jakub Zygad{\l}o}
\date{March 31, 2009}
\address{Institute of Mathematics\\
Jagiellonian University\\
Łojasiewicza 6\\
30-348 Krak\'ow, Poland.}
\email{jakub.zygadlo@im.uj.edu.pl}
 
\begin{document}
\begin{abstract}
We prove that for every pair of positive integers $a$, $b$ there exists a number $c_0$ such that for every $c\geq c_0$ one can find a tame polynomial automorphism of $\mathbb{C}^3$ with multidegree equal to $(a,b,c)$.
\end{abstract}
\maketitle

Let $F=(F_1,\ldots,F_n)$ be a polynomial automorphism of $\CC^n$ (although the results in this paper are valid over an algebraically closed field $k$ of characteristic 0, see for ex. \cite{Es} for general information about polynomial automorphisms). We define the multidegree of $F$ as $\mdeg F:=(\deg F_1,\ldots,\deg F_n)$, where $\deg$ means the total degree of a polynomial. It is an open question whether every $n$-tuple $(d_1,\ldots,d_n)$ of positive integers is a multidegree of some polynomial automorphism of $\CC^n$. 
The case $n=1$ is obvious, since $\deg F=1$ if $F$ is a polynomial automorphism. In the case $n=2$, due to Jung-van der Kulk theorem (cf. \cite{Es} and \cite{Ku} for the original proof), we must have $\deg F_1|\deg F_2$ or $\deg F_2|\deg F_1$ for $F=(F_1,F_2)$ to be a polynomial automorphism of $\CC^2$. Provided that $d_1|d_2$ (or $d_2|d_1$) one can also easily find a polynomial automorphism $F=(F_1,F_2)$ such that $\mdeg F=(d_1,d_2)$; for example $F=(x,y+x^{d_2/d_1})\circ(x+y^{d_1},y)$. The case $n=3$, however, seems much more complicated: one cannot even say (to the author's knowledge) whether there exists a polynomial automorphism with multidegree equal to $(3,4,5)$! A recent result by M. Kara\'s (\cite{Ka}) shows that if such an automorphism exists, it is not tame (see definition below and cf. \cite{SU}).

From now on we will focus our attention on the case $n=3$ and write $(a,b,c)$ instead of $(d_1,d_2,d_3)$ for the given multidegree. Without loss of generality we may assume that $c\geq b\geq a$. We will calculate an integer $c_0$ (depending on $a$ and $b$) such that for every $c\geq c_0$ one can find a polynomial automorphism of $\CC^3$ with multidegree $(a,b,c)$. Moreover, all automorphisms we introduce are tame - they are compositions of linear and elementary maps ($F=(F_1,\ldots,F_n)$ is elementary if for some $j\leq n$ and a polynomial $g$ in $n-1$ variables $F_j=X_j+g(X_1,\ldots,\widehat{X_j},\ldots,X_n)$ and $F_i=X_i$ for $i\neq j$). First, let us prove some simple facts about multidegrees.

\begin{fa}\label{f1}
Let $c\geq b\geq a$ and
\begin{enumerate}
\item $c=ka+lb$ for some $k,l\geq 0$ or
\item $a|b$ or
\item $a\leq 2$.
\end{enumerate}
Then there exists a tame automorphism of $\CC^3$ with multidegree $(a,b,c)$.
\begin{proof} We will explicitly show the required automorphism:
\begin{enumerate}
\item One can take $F=(x,y,z+x^ky^l)\circ(x+z^a,y+z^b,z)$.
\item Let $b=da$ and take $F=(x,y+x^d,z)\circ(x+y^a,y,z+y^c)$.
\item Case $a=1$ is trivial, so let $a=2$. If $2|b$ we are done due to $(2)$. So let $b$ - odd and $c=b+m$, $m\geq 0$. If $m$ is even, then $c$ is a linear combination of $b$ and $2$; if $m$ is odd, then $c$ is even - in both cases we are done due to $(1)$.\qedhere
\end{enumerate}
\end{proof}
\end{fa}

Therefore in the following we assume that $c>b>a>2$.

\begin{fa}\label{f2}
Let $b>a>2$ and suppose that $\NWD(a,b)=1$. Then for $c\geq(a-1)(b-1)$ there exists a tame automorphism of $\CC^3$ with multidegree $(a,b,c)$.
\begin{proof}
It is a classical fact due to Sylvester that every integer $c$ greater than or equal to $c_0=(a-1)(b-1)$ can be expressed as a linear combination of $a$ and $b$ with positive integer coefficients. Now apply Fact \ref{f1}.(1).
\end{proof}
\end{fa}

Note that in this case $c_0=\NWW(a,b)-a-b+1$. A similar reasoning is not valid if $\NWD(a,b):=d>1$: since $d$ divides every linear combination of $a$ and $b$, we can only get that the automorphism with multidegree $(a,b,c)$ exists for $c$ divisible by $d$ and $c\geq\NWW(a,b)-a-b+d$. However, changing $c_0$ a little, we deal with the case $\NWD(a,b)>1$ in the following theorem.

\begin{tw}\label{t1}
Let $b>a>2$ and set $r:=\min\{b-1,(a-1)(\floor{b/a}+1)\}$. Then for $c\geq\NWW(a,b)-r$ there exists a tame automorphism of $\CC^3$ with multidegree equal to $(a,b,c)$.
\begin{proof}
{\bf Step 1.} We will first show an automorphism for $c\geq\NWW(a,b)-a$.
Let $e:=\NWW(a,b)$ and write $c=e+(k-1)a+m$ for some $k\geq 0$ and $0\leq m<a$. Let us take 
\begin{eqnarray*}
&&F_1(x,y,z):=(x+z^a+z^m,y+z^b,z)\\
&&F_2(x,y,z):=(x,y,z+x^k(x^{e/a}-y^{e/b}))
\end{eqnarray*}
and consider the composition 
$$F:=F_2\circ F_1=(x+z^a+z^m,y+z^b,z+(x+z^a+z^m)^kv(x,y,z))$$
Then
$v(x,y,z)=\sum_{i=1}^{e/a}\binom{e/a}{i}(x+z^m)^iz^{e-ia}-\sum_{j=1}^{e/b}\binom{e/b}{j}y^jz^{e-jb}$ and, provided that $m>0$, the highest order monomial in $v$ is $e/a\,z^mz^{e-a}$. Therefore for $m>0$ the degree of the third coordinate of $F$ equals $ka+m+e-a=c$. If $m=0$, then $c=(e/a+k-1)\,a$ and we can take $F=(x,y,z+x^{e/a+k-1})\circ(x+z^a,y+z^b,z)$.

{\bf Step 2.}
Thanks to step 1, we must only find an automorphism $F$ with multidegree $(a,b,c)$ for $e-r\leq c<e-a$. Let us take $m:=b+c-e$ (note that $0<m<b$) and 
\begin{eqnarray*}
&&F_1(x,y,z):=(x+z^a,y+z^b+z^m+u(x,z),z)\\
&&F_2(x,y,z):=(x,y,z+x^{e/a}-y^{e/b})
\end{eqnarray*}
where $u$ is a polynomial of degree $<b$. Consider the composition
$$F:=F_2\circ F_1=(x+z^a,y+z^b+z^m+u(x,z),z+v_1(x,z)-v_2(x,y,z))$$
where
$v_1(x,z)=\sum_{i=1}^{e/a}\binom{e/a}{i}x^iz^{e-ia}$, $v_2(x,y,z)=\sum_{j=1}^{e/b}\binom{e/b}{j}(y+z^m+u(x,z))^jz^{e-jb}$.
Now let $u(x,z)=\sum_{k=1}^{\floor{b/a}}u_kx^kz^{b-ka}$ and note that
$$u(x,z)^j=\sum_{s=j}^{j\floor{b/a}}\sum_{k_1+\ldots+k_j=s}u_{k_1}\cdot\ldots\cdot u_{k_j}x^sz^{jb-sa}$$
where all the indices $k_1,\ldots,k_j$ vary from 1 to $\floor{b/a}$. One easily checks that the sum $\sum_{j=1}^{\floor{b/a}}\binom{e/b}{j}u(x,z)^jz^{e-jb}$ contains all monomials of the form $x^iz^{e-ia}$ and degree higher than $e-b$ appearing in $v_2$. Consequently, for $i\leq\floor{b/a}$, coefficient of the term $x^iz^{e-ia}$ in $v_1-v_2$ equals
$$\binom{e/a}{i}-\sum_{j=1}^{\floor{b/a}}\binom{e/b}{j}\sum_{k_1+\ldots+k_j=i}u_{k_1}\cdot\ldots\cdot u_{k_j}$$
and becomes zero if $e/b\,u_i=\binom{e/a}{i}-\sum_{j=2}^{\floor{b/a}}\binom{e/b}{j}\sum_{k_1+\ldots+k_j=i}u_{k_1}\cdot\ldots\cdot u_{k_j}$. Because on the right hand side the indices $k_1,\ldots, k_j$ are smaller than $i$, the above formula allows a recursive definition of $u_i$, $1\leq i\leq\floor{b/a}$:
\begin{align*}
u_1&:=b/a,\\
u_2&:=b/e\,\Big(\binom{e/a}{2}-\binom{e/b}{2}u_1^2\Big),\\
&\ldots\\
u_i&:=b/e\,\Big(\binom{e/a}{i}-\sum_{j=2}^{e/b}\binom{e/b}{j}\sum_{\tiny\begin{array}{c}k_1+\ldots+k_j=i\\1\leq k_1,\ldots,k_j\leq\floor{b/a}\end{array}}u_{k_1}\cdot\ldots\cdot u_{k_j}\Big)
\end{align*}
Specifying $u_i$ this way guarantees that the degree of the third coordinate of $F$ is equal to the maximum of:
$\deg(yx^{d-b})$, $\max\{\deg(x^iz^{d-ia}): \floor{b/a}<i\leq d/b\}$ and $\max\{\deg(y+z^m)^lu(x,z)^{j-l}z^{d-jb}: 1\leq j\leq d/a, 1\leq l\leq j\}$. An easy calculation shows that this maximum simplifies to $\max\{d-b+1, d-(a-1)(\floor{b/a}+1), m+d-b\}=\max\{d-r,c\}=c$.
\end{proof}
\end{tw}

\begin{uw}\label{uw1}
If in the above theorem we get $r=b-1$, then one can find an automorphism with multidegree $(a,b,c)$ for $c\geq e-b$ (since $b|e-b$).
\end{uw}

\begin{df}
A pair $(a,b)$ of positive integers such that $a<b$, $a\nmid b$ will be called \emph{tame} iff for every $c>b$ there exists a tame automorphism of $\CC^3$ with multidegree equal to $(a,b,c)$.
\end{df}

\begin{wn}
A pair $(2n,kn)$, where $k\geq 3$ - odd is tame for $n\geq\frac{k-1}{2}$.
\begin{proof}
Due to Fact \ref{f1}.(3) we can assume $n\geq 2$. We apply the Theorem and get the existence of tame automorphism with multidegree $(2n,kn,c)$, provided that $c\geq 2kn-r$, where $r=\min\{kn-1, (2n-1)(\floor{kn/2n}+1)\}=\min\{kn-1, (2n-1)(k+1)/2\}=\min\{kn-1,kn+n-\frac{k+1}{2}\}$. Following Remark \ref{uw1}, we are done if $n-\frac{k+1}{2}\geq-1$ \ie $n\geq\frac{k-1}{2}$.
\end{proof}
\end{wn}

\begin{wn}
Pairs $(2n,3n)$ and $(2n,5n)$ are tame for all $n\geq 1$. 
\end{wn}

\begin{uw}
In the case $n=\frac{k-3}{2}$ there are only following triples "left": $(4,14,15)$, $(6,27,28)$, $(8,44,45)$, ..., $(2n,(2n+3)n,(2n+3)n+1)$.
\end{uw}

\end{document}